\documentclass{article}
\usepackage{amsmath,amssymb}
\usepackage{euscript}
\usepackage{latexsym}
\usepackage[pdftex]{color}
\usepackage[matrix,arrow,curve]{xy}
\usepackage[english]{babel}
\usepackage{wrapfig}
\usepackage{mathrsfs}

\newtheorem{thm}{\bf Theorem}[section]

\newtheorem{prop}[thm]{\bf Proposition}
\newtheorem{corollary}[thm]{\bf Corollary}
\newtheorem{definition}{\bf Definition}[section]
\newtheorem{example}[definition]{\bf Example}

\newcounter{abcd}[section]
\addtocounter{abcd}{1}

\newcommand{\?}{\mathstrut}

\def\le{\leqslant}
\def\ge{\geqslant}

\begin{document}

\begin{center}
{\large
	COHOMOLOGY RINGS OF PRECUBICAL SETS
}
\\
\medskip
{
Lopatkin V.E.
}\\
\end{center}

\begin{abstract}
The aim of this paper is to define the structure of a ring on a graded cohomology group of a precubical set in coefficients in a ring with unit.
\end{abstract}

\section*{Introduction}
Let $G$ be the homologous system Abelian groups over a precubical set $X$ [1], then for any integral $n \ge 0$, $H_n(X;G)$ are defined by values of satellites of the colimit functor ${\varprojlim}^n: \mathrm{Ab}^{\left(\square_+/X\right)^{op}} \to \mathrm{Ab}$, here $\square_+/X$ is a category of singular cubes of a precubical set $X$, $\mathrm{Ab}$ is the category of Abelian groups and homomorphisms, further for any small category $\mathscr{C}$ we denote by $\mathscr{C}^{op}$ the opposite category and finally $\mathrm{Ab}^{\left(\square_+/X\right)^{op}}$ is the category of functors from $\left(\square_+/X\right)^{op}$ to $\mathrm{Ab}$. This observation is generalizing the Serre's spectral sequence for precubical sets [1]. For the cohomology groups there exist a opposite statement.
\par A cohomologous system over a precubical set we define as a functor on a category of singular cubes. In general, values of this functor on morphisms are not isomorphisms.
\par Suppose that the cohomologous system take constant values which are any ring $R$ then we can to define a structure of a ring over a graded cohomology group with coefficients in this system.  
\par The aim of this paper is to define the structure of a graded ring over a graded cohomology group of precubical sets with coefficients in the cohomologous system witch is taken a constant value. The basic result of this paper is Theorem 4.4.
\par We use following notations. The category of sets and maps we denote by $\mathrm{Ens}$, $\mathrm{Ab}$ is the category of Abelian groups and homomorphisms and $\mathrm{Ring}$ is the category of rings and ring's homomorphisms which are save the unit.

\section{Precubical Sets}
\begin{definition}
A precubical set $X = (X_n, \partial_i^{n, \varepsilon})$ is a sequence of sets $(X_n)_{n \in \mathbb{N}}$ with a fimile of maps $\partial_i^{n,\varepsilon}: X_n \to X_{n-1}$, defined for $i \le i \le n, \varepsilon \in \{0,1\}$, for which the following diagrams is commutative for all $\alpha, \beta \in \{0,1\},$ $n \ge 2$, $1 \le i < j \le n:$
$$
  \xymatrix{
Q_n \ar@{->}[r]^{\partial^{n,\beta}_j} \ar@{->}[d]_{\partial^{n,\alpha}_i}  & Q_{n-1} \ar@{->}[d]^{\partial^{n-1,\alpha}_i} \\
Q_{n-1} \ar@{->}[r]_{\partial^{n-1,\beta}_{j-1}} & Q_{n-2}
}
$$ 
\end{definition}
\par Let $\square_+$ be a category consisting of finite sets $\mathbb{I}^n = \{0,1\}$ ordered as the Cartesian power of $\mathbb{I}$. Any morphism of the $\square_+$ is defined as an ascending map which admits a decomposition of the form $V^{k,\varepsilon}_i: \mathbb{I}^{k-1} \to \mathbb{I}^k$ where
$$V^{k,\varepsilon}_i (u_1, \ldots, u_{k-1}) = (u_1, \ldots, u_{i-1},\varepsilon, u_i, \ldots, u_{k-1}), \qquad \varepsilon \in \{0,1\}, \quad 0 \le i \le k.$$
here $\varepsilon \in \{0,1\}$, $1 \le i \le k$. Also we'll denote maps $V^{n,\varepsilon}_i$ by $V^\varepsilon_i$. 
\par It well know [1] that any precubical set $X$ is a functor $X:\square_+^{op} \to \mathrm{Ens}$.
\par Let $H$ be a ordered subset $\{h_1,\ldots, h_p\}$ of the set $\{1,2,\ldots, n\}$. Let us define a map $\lambda_H^\varepsilon: \mathbb{I}^p \to \mathbb{I}^n$ by the following formula
$$\lambda_h^\varepsilon \left( u_1, \ldots, u_p \right)  =  \left( v_1, \ldots, v_n \right),$$
where $v_i = \varepsilon$, if $i \notin H$, and $v_{h_r} = u_r$, $r = 1, \ldots, p$.
\begin{prop}
Suppose that we have a subset $H = \{h_1, \ldots, h_p\}$ of the set $\{1,2,\ldots, n\}$. Let us define following sets; $\widehat{H}_\mu = \{h_1, \ldots,h_{\mu-1},h_{\mu+1},\ldots, h_p\}$, $\widetilde{H}_\mu = \{h_1,\ldots,h_{\mu-1}, h_{\mu+1}-1,\ldots,h_p-1\}$. Further, let $H_j$ be a $\{h_1, \ldots, h_r, h_{r+1} -1, \ldots, h_p-1\}$ if $j \notin H$ and $h_r < j< h_{r+1}$. There are following formulas for $\varepsilon, \eta \in \{0,1\}$
$$\lambda_H^\eta \circ  V_\mu^\varepsilon  =  V_{h_\mu}^\varepsilon \circ \lambda_{\widetilde{H}_\mu}^\eta;$$
$$\lambda_H^\varepsilon \circ  V_\mu^\varepsilon  = \lambda_{\widehat{H}_\mu}^\varepsilon; $$
$$\lambda_H^\varepsilon =  V_j^\varepsilon \circ \lambda_{H_j}^\varepsilon. $$
In this case, $H_j$ and $\widetilde{H}_\mu$ are subsets of set $\{1,2, \ldots, n-1\}$.
\end{prop}
\textbf{Proof} This proposition was proved in [2, Proposition 9.3.4]

\section{A Diagonal Inclusion}
In this section we'll introduce a diagonoal inclusion and show that this inclusion is the chain map. 
\par First let us intruduce some notices from [1].
\par Let $X = (X_n, \partial^{n,\varepsilon}_i)$ be the presubical set, let $\square_+[X_p] = L(X_p)$ for $p \ge 0$ be free Abelian group and $\square_+[X_p] = 0$ for $p < 0$. Assume that $D^\varepsilon_i = L(\partial_i^\varepsilon): \square_+[X_p] \to \square_+[X_{p-1}]$. Further let us define homomorphisms
$$D : \square_+[X_p] \to \square_+[X_{p-1}], \quad p \ge 1,$$
by the formula
$$D = \sum\limits_{i=1}^p (-1)^i \left(D_i^1 - D_i^0 \right).$$
\par Let us assume that $\square_+ [X] = \bigoplus\limits_{p \ge 0} \square_+ [X_p]$ be the direct sum of groups $\square_+ [X_p]$.

\par Following [1], identify cubes $f \in X_p$ with corresponding natural transformations $\widetilde{f}:h_{\mathbb{I}^p} \to X$ which are called \emph{singular cubes}. Thus, singular $p$-cubes are elements of the group $\square_+[X_p]$.
\par Let us consider functor morphisms $h_{\lambda_H^\varepsilon}: h_{\mathbb{I}^p} \to h_{\mathbb{I}^n}$, it's hard to see that the homomorphism $D$ can define by the following corresponding $D^\varepsilon: f \mapsto f \circ h_{V^\varepsilon_i}$. It is clear that the $f \circ h_{V^\varepsilon_i}$ is define any face of the singular $p$-cube. There are rules of commutation functor morphisms $h_{\lambda_H^\varepsilon}$ with the homomorphism $D$ in the following proposition which is a modification of proposition 1.1
\begin{prop}
Let us assume that we have a ordered subset $G = \{g_1, \ldots, g_p\}$ of set $\{1,2,\ldots, n\}$. Suppose that $\widehat{G}_\mu= \{g_1, \ldots,g_{\mu-1},g_{\mu+1},\ldots, g_p\}$ and $\widetilde{G}_\mu = \{g_1,\ldots,g_{\mu-1}, g_{\mu+1}-1,\ldots,g_p-1\}$. Furhter, suppose that $G_j= \{g_1, \ldots, g_r, g_{r+1} -1, \ldots, g_p-1\}$ if $j \notin G$ and $g_r < j < g_{r+1}$. Let us assume that we have a precubical set $X = \left(X_n, \partial_i^{n,\varepsilon}\right)$, let $f:h_{\mathbb{I}^n} \to X$ be a singular $n$-cube. There are following formulas for $\varepsilon, \eta \in \{0,1\}:$
$$D_{\mu}^{\varepsilon}\left(f\circ h_{\lambda_{G}^{\eta}}\right)= 
D_{g_{\mu}}^{\varepsilon}\left(f\right)\circ h_{\lambda_{\tilde{G}_{\mu}}^{\eta}}  \eqno(1)$$

$$D_{\mu}^{\varepsilon}\left(f\circ h_{\lambda_{G}^{\varepsilon}}\right)  =  f \circ h_{\lambda_{\widehat{G}_\mu}^\varepsilon}  \eqno(2)$$

$$f \circ h_{\lambda_G^\varepsilon}  =  D_j^{\varepsilon}(f)\circ h_{\lambda_{G_j}^{\varepsilon}} \eqno(3)$$
We assumed that $G_j$ and $\widetilde{G}_\mu$ are ordered subsets of set $\{1,2,\ldots,n-1\}$.
\end{prop}
\textbf{Proof.} From proposition 1.1 it follows that there are following formulas
$$h_{\lambda_G^\eta} \circ h_{V_\mu^\varepsilon} = h_{V_{g_\mu}^\varepsilon} \circ h_{\lambda_{\widetilde{G}_\mu}^\eta} ;$$
$$ h_{\lambda_G^\varepsilon} \circ h_{V_\mu^\varepsilon} = h_{\lambda_{\widehat{G}_\mu}^\varepsilon}; $$
$$h_{\lambda_G^\varepsilon} =h_{\lambda_{G_j^\varepsilon}} \circ h_{V^\varepsilon_j} . $$
\par Multiplying both sides by $f$, we complete the proof (see the commutative diagramm). 
$$
  
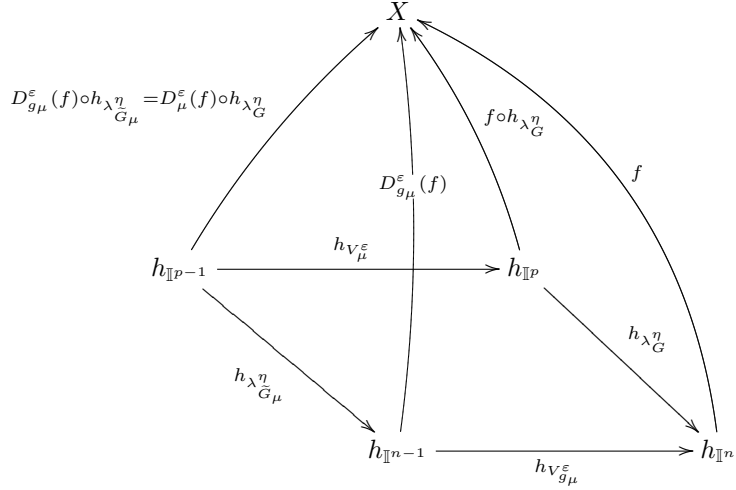
\xymatrix{
&& X &&&\\
&&&&&&&\\
&&&&&&&\\
h_{\mathbb{I}^{p-1}} \ar@{->}[rrdd]_{h_{\lambda^\eta_{\widetilde{G}_\mu}}}   \ar@/^/@{->}[uuurr]^{D^\varepsilon_{g_\mu}(f) \circ h_{\lambda_{\widetilde{G}_\mu}^\eta} = D^\varepsilon_\mu(f) \circ h_{\lambda_G^\eta}} \ar@{->}[rrr]^{h_{V^\varepsilon_\mu}} &&& h_{\mathbb{I}^{p}} \ar@{->}[rrdd]^{h_{\lambda^\eta_G}} \ar@/_/@{->}[uuul]_{f \circ h_{\lambda^\eta_G}} && \\
&&&&&&&\\
&& h_{\mathbb{I}^{n-1}} \ar@/_/@{->}[uuuuu]|(.6){D^\varepsilon_{g_\mu}(f)} \ar@{->}[rrr]_{h_{V^\varepsilon_{g_\mu}}} &&& h_{\mathbb{I}^{n}} \ar@/_2.5pc/@{->}[uuuuulll]_{f}
}
$$

\par It well know (see [3]) that the tensor product $\square_+ [X] \otimes \square_+ [X]$ of the chain complex $\square_+ [X]$ with itself is the chain complex $\square_+ [X \otimes X]$, where
$$\square_+ [(X \otimes X)_n] = \bigoplus\limits_{p+q=n} \square_+ [X_p] \otimes \square_+ [X_q],\eqno(4)$$
and bound operators is defined over generators $x\otimes x'$ by the formula
$$\partial (x \otimes x') = \partial x \otimes x' + (-1)^{\mathrm{dim}\,x}x \otimes \partial x'.\eqno(5)$$\\
\begin{prop}
Let $X \in \square_+^{op}\mathrm{Ens}$ be a precubical set, let us assume that $\square_+[X]$ be aforesaid chain complex. Further let $\square_+ [X \otimes X]$ be the tensor product of the chain complex $\square_+ [X]$ with itself which defined by the formulas (4), (5). A map $\Delta$ (diagonal inclusion) which defined by the formula for any singular cub $f: h_{\mathbb{I}^n} \to X$: 
$$\Delta (f) = \sum\limits_G \varrho_{GK} \left(f \circ h_{\lambda_G^0} \right) \otimes \left( f \circ h_{\lambda_K^1} \right) , $$
is the chain map. Here $K$ is the complement of a set $G = \{g_1,\ldots, g_p\} \subseteq \{1,2,\ldots, n\}$, $varrho_{GK}$ is a signature of a permutation $GK$  of integral numbers $1,2,\ldots, n$. The summation is taken over all ordered subsets $G$ of set $\{1,2,\ldots n\}$.
\end{prop}
\textbf{Proof} This poropsition was proved in [2, Proposition 9.3.5]

\section{Cohomology of Precubical Sets with Coefficients in a Cohomologous System of Rings}
\begin{definition}
A cohomologous system of rings and a cohomologous system of Abelian groups over a precubical set $X \in \square_+^{op}\mathrm{Ens}$ are some functors $\mathscr{R}: \square_+/X \to \mathrm{Ring}$ and $\mathscr{G}: \square_+/X \to \mathrm{Ab}$, respectively.
\end{definition}
\par Let us consider Abelian grpups $^n\square_+[X,\mathscr{G}] = \prod\limits_{\vartheta \in X_n} \mathscr{G} (\vartheta)$. Let us define differentials $\delta^{n,\varepsilon}_i: \?^n\square_+[X,\mathscr{G}] \to \?^{n+1}\square_+[X,\mathscr{G}]$ as homomorphisms making following diagrams commutative
$$
 \xymatrix{
  \prod\limits_{\vartheta \in X_n} \mathscr{G}( \vartheta) \ar@{->}[d]_{\mathrm{pr}_{\vartheta \circ V^{n+1,\varepsilon}_i}} \ar@{->}[rrr]^{\delta^{n,\varepsilon}_i} &&& \prod\limits_{\vartheta \in X_{n+1} }\mathscr{G}(\vartheta) \ar@{->}[d]^{\mathrm{pr}_\vartheta}  \\
\mathscr{G} \left(\vartheta \circ V^{n+1,\varepsilon}_i \right) \ar@{->}[rrr]_(.55){\mathscr{G}\left(V^{n+1,\varepsilon}_i: \vartheta V^{n+1,\varepsilon}_i \to \vartheta \right)}  &&& \mathscr{G}(\vartheta)
}
$$
\begin{definition}
Let $X$ be a precubical set, let $\mathscr{G}: \left.\square_+  \right/ X \to \mathrm{Ab}$ be a cohomologous system of Abelian groups over this precubical set $X$. A cohomology groups $H^n(X;\mathscr{G})$ of this precubical set $X$ with coefficients in $\mathscr{G}$ are $n$-th cohomology groups of a chain complex $^*\square_+[X,\mathscr{G}]$ consisting of abelian groups 
$$^{n}\square_+[X,\mathscr{G}] = \prod\limits_{\sigma  \in X_n} \mathscr{G}(\sigma )$$
and differentials 
$$\delta^n = \sum\limits_{i=1}^{n+1}(-1)^i (\delta^{n,1}_i- \delta^{n,0}_i).$$
\end{definition}
\par Suppose that the cohomologous system of rings $\mathscr{R}: \square_+/X \to \mathrm{Ring}$ over a precubical set $X$ take a constant value which is a ring $R$ with a unity. Considering an additive component of the ring $R$ we can examine a cohomology groups $H^*(X;R)$ with coefficient in the ring $R$.
\par Let $\?^* \square_+[X;R]$ be a cochain complex. Following [2, \S 5.7, 5.7.27] let us consider the homomorphism
$$\pi: \?^*\square_+ [X; R] \otimes_{R} \?^*\square_+ [X; R] \to \?^*\square_+ [X \otimes X; R],$$
which defined by the formula
$$\left( \pi (u \otimes u')  \right) (c \otimes c')  = \eta \left(u(c) \otimes_{R} u'(c')  \right),$$
here $c,c' \in \square_+ [X]$, $u,u' \in \?^*\square_+ [X; R ]$ and $\eta: R\otimes_{R} R \to R$ is an isomorphism of rings wich defined by the following formula
$$\eta \left(u(c) \otimes u'(c') \right) = u(c) \cdot u'(c'),$$
this prodoct is the multiplication operation in the ring $R$. 
\par From [2, Proposition 5.7.28] follow that the homorphism $\pi$ is the cochain map. Thus it's not hard to see that a map
$$\smile = \Delta^* \pi : \?^*\square_+[X;R] \otimes_{R} \?^*\square_+[X;R] \to \?^*\square_+[X;R] $$
is the cochain map because from proposition 2.2 follows that the map $\Delta^*$ is the cochain map. It means that the $\smile$ generate some a product in $H^*(X;R)$. Thus we have the following
\begin{thm}
The graded group $H^*(X;R)$ with afore--mentioned $\smile$-product is a ring.
\end{thm}
\par Let us describe the $\smile$-product over cochains. Let $\varphi \in \?^p \square_+ [X; R]$ and $\psi \in \?^q \square_+ [X; R]$ are cochains. Let $u  \in X_{p+q}$ be a $p+q$-cube. We have a formula
$$(\varphi \smile \psi) (u) = \sum\limits_G \varrho_{GK} \varphi\left(u \circ h_{\lambda_G^0} \right) \cdot \psi\left(u \circ h_{\lambda_K^1} \right),\eqno(6)$$
Here $G = \{g_1,\ldots, g_p\} \subseteq \{1,2,\ldots, n\}$, $\varrho_{GK}$ is a signature of a permutation $GK$  of integral numbers $1,2,\ldots, n$. The summation is taken over all ordered subsets $G$ of set $\{1,2,\ldots n\}$.
\par The notices of form $u \circ h_{\lambda_G^0}$ we also denote by $uh_{\lambda_G^0}$.

\section{Properties of the Precubical Cohomology Ring}
Here we will enumerate and we'll proof algebraic properties of the $\smile$-product in the ring $H^*(X;R)$.
\begin{thm}
The $\smile$-product of cochains in the ring $\?^*\square_+[X;R]$ is associative and distributive with respect to the addition. If the ring $R$ has left (right, two--sided) unit then the ring $\?^*\square_+[X;R]$ has same unit.
\end{thm}
\textbf{Proof.} From associative and distributive of the product in the ring $R$ follows associative and distributive of product in the ring $\?^*\square_+[X;R]$. Further, let $1$ --- be a left unit of the ring $R$ and let $\iota$ be a cochain which take each the $0$-cube of the $\square_+[X]$ to $1$. It's not hard to see that for any cochain $\xi$ there is the following equality $\iota \smile \xi = \xi$. In the same way we'll get the proof of this theorem if $1$ is right or two--sided unit.\\

\par The cochain complex $\?^*\square_+[X;R]$ with $\smile$-product is a graded ring.

\begin{thm}
For $\varphi \in \?^p\square_+[X;R]$ и $\psi \in \?^q\square_+[X;R]$ there is the following formula
$$\delta \left( \varphi \smile \psi \right) = \delta \varphi \smile \psi + (-1)^p \varphi \smile \delta \psi $$
\end{thm}
\textbf{Proof.} We have
$$\left( \delta \varphi \smile \psi  \right) (f) = \sum\limits_{G}\varrho_{GK} \, \left( \delta \varphi \right) \left(f h_{\lambda_G^0} \right) \cdot \psi \left(f h_{\lambda_K^1} \right)=$$
$$ = \sum\limits_{G} \varrho_{GK} \left( \sum\limits_{\mu =1}^{p+1} (-1)^\mu \left( \left( \delta_\mu^1 \varphi \right) \left(f h_{\lambda_G^0} \right) - \left( \delta_\mu^0 \varphi \right) \left(f h_{\lambda_G^0} \right) \right) \right) \cdot \psi \left(f h_{\lambda_K^1}  \right),$$

$$\left( \varphi \smile \delta \psi  \right) (f) = \sum\limits_{G}\varrho_{GK} \, \left( \varphi \right) \left( f h_{\lambda_G^0} \right) \cdot \left( \delta  \psi \right) \left( f h_{\lambda_K^1} \right)=$$
$$ = \sum\limits_{G} \varrho_{GK} \, \varphi \left(f h_{\lambda_G^0} \right) \cdot \left( \sum\limits_{\eta =p}^{p+q+1} (-1)^\eta \left( \left( \delta_\eta^1 \psi \right) \left( f h_{\lambda_K^1} \right) - \left( \delta_\eta^0 \psi \right) \left( f h_{\lambda_K^1} \right) \right) \right),$$
here $G \subset \{1,2, \ldots, p+q+1 \}$, $G = (h_1, \ldots, h_{p+1})$ and $K$ is the complement of the set $G$.
\par From the diagram
$$
  \xymatrix{
h_{\mathbb{I}^{p+q+1}} \ar@{->}[dd]_f && \ar@{->}[ll]_{h_{\lambda_G^\xi}}  h_{\mathbb{I}^{p+1}}  \ar@{->}[ddll]| { f h_{\lambda_G^\xi} } \\
&&&&\\
X && h_{\mathbb{I}^p} \ar@{->}[ll]^{D^\varepsilon_\mu  \left(f h_{\lambda_G^\xi} \right)} \ar@{->}[uu]_{h_{V^\varepsilon_\mu}}
}
$$
and proposition 2.1 it follows that, we have
$$\left( \delta \varphi \smile \psi  \right) (f) = \sum\limits_{G}\varrho_{GK} \left( \sum\limits_{\mu=1}^{p+1} (-1)^\mu \left( \varphi \left[D^1_{h_\mu} \left(f h_{\lambda_{\widetilde{G}_\mu}^0} \right) \right]   \right)  -  \varphi \left[f h_{\lambda_{\widehat{G}_\mu}^0} \right]   \right)  \cdot \psi \left[ f h_{\lambda_K^1}  \right],$$

$$\left( \varphi \smile \delta \psi  \right) (f) = \sum\limits_{G}\varrho_{GK} \, \varphi \left[ f h_{\lambda_G^0 }  \right] \cdot \left( \sum\limits_{\eta=p}^{p+q+1} (-1)^\eta \left( \psi \left[f   \lambda_{\widehat{K}_\eta}^1  \right]   \right)  -  \psi \left[ D_{k_\eta}^0   \left( f h_{\lambda_{\widetilde{K}_\eta}^1} \right) \right]   \right).$$

\par Let $\check{K}_\mu$ be a complement of the set $\widehat{G}_\mu$. Let us consider a sum $\left( \delta \varphi \smile \psi  \right) (f) + (-1)^p \left( \varphi \smile \delta \psi  \right) (f)$. It's not hard to see that $\varphi \left[f h_{\lambda_{\widehat{G}_\mu}^0 } \right] \cdot \psi \left[ f h_{\lambda_K^1} \right]$ will appear twice; in the first place it will appear as a result of a deletion the $g_\mu$ from the $G$ in the component $(G,K)$ and in the second place it will appear as a result of a deletion the $g_\mu$ from the $\check{K}_\mu$ in the component $(\widehat{G}_\mu, \check{K}_\mu)$. In the first place $\varphi \left[f h_{\lambda_{\widehat{G}_\mu}^0} \right] \cdot \psi \left[ f h_{\lambda_K^1} \right]$ hase a sign $\varrho_{GK} (-1)^{\mu+1}$, , further, in the second place it hase a sign $\varrho_{\widehat{G}_\mu \check{K}_\mu} (-1)^p (-1)^\alpha$, here $k_\alpha < g_\mu < k_{\alpha +1}$. But we have
$$\varrho_{\widehat{G}_\mu \check{K}_\mu} = (-1)^{p-\mu+\alpha} \varrho_{GK},$$
it means that the $\varphi \left[f h_{\lambda_{\widehat{G}_\mu}^0} \right] \cdot \psi \left[ f h_{\lambda_K^1 } \right]$ will appear twice with different signs. So that we have
$$\left( \delta \varphi \smile \psi  \right) (f) + (-1)^p \left( \varphi \smile \delta \psi  \right) (f) = $$
$$= \sum\limits_{G}\varrho_{GK} \left( \sum\limits_{\mu=1}^{p+1} (-1)^\mu \varphi \left[D^1_{g_\mu}   \left( f  h_{\lambda_{\widetilde{G}_\mu}^0} \right) \right]   \cdot \psi \left[ f h_{\lambda_K^1}  \right]  + \right.$$
$$ \left. \qquad + (-1)^p  \sum\limits_{\eta=p}^{p+q+1} (-1)^{\eta+1} \varphi \left[ f h_{\lambda_G^0} \right] \cdot \psi \left[D^0_{k_\eta}   \left( f h_{\lambda_{\widetilde{K}_\eta}^1} \right) \right]  \right). \eqno(7)$$

\par From other side we have
$$\left(\delta (\varphi \smile \psi)  \right) (f) = \sum\limits_{i=1}^{p+q+1} (-1)^i \left( (\varphi \smile \psi ) \left( D_i^1 f \right) - (\varphi \smile \psi ) \left( D_i^0 f \right)  \right) = $$
$$=  \sum\limits_{i=1}^{p+q+1} (-1)^i \sum\limits_F \varrho_{FT}\, \left( \varphi \left[ D_i^1(f) h_{\lambda_F^0} \right] \cdot \psi \left[ D^1_i(f) h_{\lambda_T^1} \right] - \right.$$    
$$ \left. -\varphi \left[ D_i^0 (f) h_{\lambda_F^0}   \right] \cdot \psi \left[ D_i^0 (f) h_{\lambda_T^1}  \right]\right), \eqno(8)$$
here $F$ is an ordered subset of the set $\{1,2,\ldots, p+q \}$ and $T$ is its complement.
\par Using (1) -- (3) of proposition 2.2, and assume that 
$$F = \begin{cases}\widetilde{G}_j; \, \mbox{если} \quad j \in G\\ G_j;\, \mbox{если}  \quad j \notin G   \end{cases} \qquad \mbox{и} \qquad T =  \begin{cases}{K}_j; \, \mbox{если} \quad j \in G\\ \widetilde{K}_j;\, \mbox{если}  \quad j \notin G   \end{cases}$$
we get a bijection between triples $(F,T,i)$ and $(G,K,j)$ here $i=j$. It means that we have a bijection between (7) and (8) up to the sign.
Let us prove that this signs are equal. We must check the following equation
$$ (-1)^\mu \varrho_{GK} = (-1)^{h_\mu} \varrho_{\widetilde{G}_\mu K_\mu}, \qquad (-1)^\eta \varrho_{GK} = (-1)^{k_\eta} \varrho_{{G}_\mu \widetilde{K}_\mu}.$$
Let us compare followings permutations
$$GK: \quad g_1, \ldots, g_{\mu-1}, g_\mu, \ldots, g_p, k_1, \ldots, k_\alpha, k_{\alpha+1}, \ldots, k_q$$
and
$$\widetilde{G}_\mu K_{h_\mu}: \quad g_1, \ldots, g_{\mu-1}, g_{\mu+1}-1, \ldots, g_p-1, k_1, \ldots, k_\alpha, k_{\alpha+1}-1, \ldots, k_q-1,n.$$
It's not hard to see that following permutations
$$g_\mu,\ldots, g_p, k_{\alpha+1}, \ldots, k_q \quad \mbox{and} \quad g_{\mu+1}-1, \ldots, g_p-1, k_{\alpha+1}-1, \ldots, k_q -1,n$$
have same signs, because we can get from first to second permutation by two steps: in the first step, we add $1$ to all numbers, so we get $g_{\mu+1}, \ldots, g_p, k_{\alpha+1}, \ldots, k_q,g_\mu$, and in the second step we transfer $g_\mu$ in the beginning. Each of this steps multiply the sing by $(-1)^{n-g_\mu}$. It means that sings of the last permutation are differents with respect to the $(-1)^\alpha$. Here $\alpha$ is a number of $k$ which are smaller than $g_\mu$, so that $\alpha = g_\mu - \mu$ and we complete to proof the first equation. In just the same way we can to proof the second equation.
\begin{flushright}
Q.E.D.
\end{flushright}

\par Let a cochain complex is a graded ring with respect to any product, then this cochain complex is said [2] to be a \emph{cochain ring}, if this product satisfy theorem 4.2. From theorem 4.2 we get the following
\begin{corollary}
If $\varphi$ and $\psi$ are cocycles, then $\varphi \smile \psi$ is a cocycle. Moreover if $\xi$ is a coboundary and $\zeta$ is a cocycle then $\xi \smile \zeta$ is a coboundary.
\end{corollary}  
\textbf{Proof.} Indeed, using theorem 4.2, we get
$$\delta (\varphi \smile \psi) = \delta (\varphi) \smile \psi + (-1)^{\mathrm{dim} \varphi} \varphi \smile (\delta \psi) = 0+ 0=0.$$
\par Let us suppose that $\xi = \delta \vartheta$ and let $\xi$ be a coboundary, further let $\zeta$ be a cocycle, then 
$$\delta (\vartheta \smile \zeta) = (\delta \vartheta) \smile \zeta + (-1)^{\mathrm{dim} \vartheta} \vartheta \smile (\delta \zeta) = \xi \smile \zeta.$$ 
This completes the proof of this Corollary.\\
\par Now we formulate the basic result of this paper.
\begin{thm}
A set $Z\left(X; R  \right)$  of cocycles is a subring of the ring $\?^*\square_+[X;R]$; a set $B\left( X ; R  \right)$ of coboundaries is a two--sided ideal in the ring $Z\left(X; R  \right)$. The cohomology ring $H^*(X;R)$ of the a precubical set $X \in \square_+^{op}\mathrm{Ens}$ is isomorphic to the quotient--ring $Z(X;R)/B(X;R)$. The ring $H^*(X;R)$ is a graded ring. If the ring $R$ has left (right, two--sided) unity, then the ring $H^*(X;R)$ has the same unity. 
\end{thm}
\textbf{Proof.} From Corollary 4.3 it follows that a set $Z(X;R)$ is a subring of the ring $\?^* \square_+[X;R]$ and a set $B\left( X ; R  \right)$ is a two--sided ideal in the ring $Z\left(X; R  \right)$. Further, from Definition 3.2 we get a additive isomorphism $H^*(X;R) \cong Z(X;R) / B(X;R)$. Suppose that $f,g \in H^*(X;R)$, let us consider their representatives $[f]$ and $[g]$ in $Z(X;R)$, respectively. It's not hard to see that using (6), we have that a representative of $f \smile g$ be $[f \smile g]$. It's evident that the above--cited cochain $\iota$ is a cocycle, this completes the proof of this Theorem.\\
\par Let us show that there is the following
\begin{thm}
If the ring $R$ is a commutative then the ring $H^*(X;R)$ is an anticommutative.
\end{thm} 
\textbf{Proof.} Since for any permutation $GK$ of integral numbers $1,2,\ldots, n$ there is the following equation $\varrho_{GK} = \varrho_{KG}$ then we get for any $\varphi \in \?^p\square_+[X;R]$, $\psi \in \?^q\square_+[X;R]$ the following equation
$$\varphi \smile \psi = (-1)^{pq}\psi \smile \varphi.$$
\begin{flushright}
Q.E.D
\end{flushright}

\begin{example}
Let us to calculate the cohomology ring of the torus $\mathbb{T}^2$. We present the torus $\mathbb{T}^2$ as a precubical set $\mathbb{T}^2 = \left( Q_n\mathbb{T}^2; \partial_i^{n,\varepsilon} \right)$, see the figure 1.
\par \begin{figure}[h!]
\begin{center}
\begin{picture}(10,10)
\linethickness{0,5mm}
\put(-40,0){\begin{picture}(10,10)
\put(0,0){\line(1,0){100}}
\put(100,0){\line(0,-1){100}}
\put(100,-100){\line(-1,0){100}}
\put(0,-100){\line(0,1){100}}
\put(0,0){\circle*{5}}
\put(-10,2){$A$}
\put(100,0){\circle*{5}}
\put(102,2){$B$}
\put(100,-100){\circle*{5}}
\put(102,-110){$C$}
\put(0,-100){\circle*{5}}
\put(-10,-110){$D$}
\put(60,0){\vector(1,0){5}}
\put(100,-60){\vector(0,-1){5}}
\put(40,-100){\vector(-1,0){5}}
\put(0,-40){\vector(0,1){5}}
\multiput(50,20)(0,-2){70}{\circle*{1}}
\put(51,-123){$\alpha$}
\multiput(-20,-50)(2,0){70}{\circle*{1}}
\put(-23,-47){$\beta$}

\end{picture}
}
\end{picture}
\end{center}
\vspace{28ex}
\caption{Here is shown the expanding of the torus; $DA$ is identified with $CB$ and $AB$ is identified with $DC$.}
\end{figure}
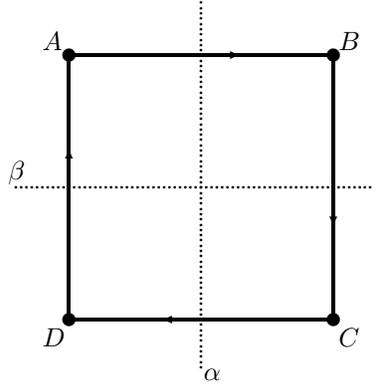 
\par So we have $Q_0\mathbb{T}^2 = \{o = A = B = C = D\}$, $Q_1 \mathbb{T}^2 = \{t_1 = DA = CB, t_2 = AB = DC\}$, $Q_2 \mathbb{T}^2 = \{\vartheta = ABCD\}$. In the figure 2 are shown values which are taken bound differentials on the one and the two--dimension cubes.
\par We have the following cochain complex
$$0 \to \mathbb{Z}^{1} \xrightarrow{\delta^0} \mathbb{Z}^{2} \xrightarrow{\delta^1} \mathbb{Z}^1 \xrightarrow{\delta^2} 0$$
\par Let us to assign $k$--dimension cochain $\vartheta^*$ to each $k$-cube $\vartheta$ of the precubical torus. This cochain $\vartheta^*$ is taken $1$ on the cube $\vartheta$ and it is taken $0$ on others cubes. We'll consider cochains which are the sum of cochains of form $\vartheta^*$. 
\par Since the following diagram is commutative
$$
\xymatrix{
&&& f \in \prod\limits_{\vartheta \in Q_n \mathbb{T}^2} \mathbb{Z} \ar@{->}[dd]_{\mathrm{pr}_{\vartheta \circ V^{n+1, \varepsilon}_i}} \ar@{->}[rrr]^{\delta^{n,\varepsilon}_i} \ar@/^1.5pc/@{->}[rrrdd]|{\left(\delta^{n,\varepsilon}_i f \right)(\vartheta)}  \ar@/_1.5pc/@{->}[rrrdd]|{f \left(\vartheta  V^{n+1,\varepsilon}_i \right)}  &&&  \prod\limits_{\vartheta \in Q_{n+1} \mathbb{T}^2} \mathbb{Z}  \ar@{->}[dd]^{\mathrm{pr}_\vartheta} \\
&&&& &&& & \\
&&& \mathbb{Z} \ar@{->}[rrr]_{\mathbb{Z}\left( V^{n+1,\varepsilon}_i: \vartheta \circ V^{n+1, \varepsilon}_i  \to  \vartheta  \right)} &&& \mathbb{Z} 
}
$$
then there exist the following equation
$$\left(\delta^{n,\varepsilon}_i f \right)(\vartheta) = f \left(\vartheta  V^{n+1,\varepsilon}_i \right).$$
From this equation it's not hard to see that the one--dimension cochain $f$ is taken different sign values on two edges of the bound of the $2$-cube (according to the sign of the orientation of this $2$-sube) then $f$ is the cocycle.(see fig. \ref{cub}).
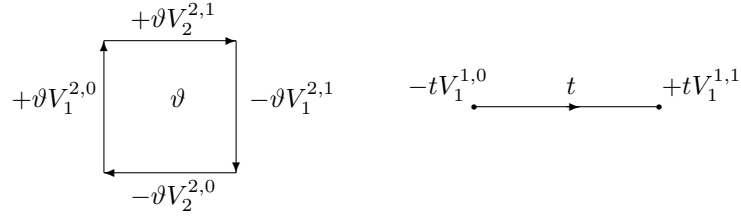
\begin{figure}[h!]
\begin{center}
\begin{picture}(15,15)
\put(-70,0){
\begin{picture}(1,1)
\put(0,0){\vector(1,0){50}}
\put(50,0){\vector(0,-1){50}}
\put(50,-50){\vector(-1,0){50}}
\put(0,-50){\vector(0,1){50}}
\put(-35,-25){$+\vartheta V^{2,0}_1$}
\put(55,-25){$- \vartheta V_1^{2,1}$}
\put(10,5){$+\vartheta V^{2,1}_2 $}
\put(10,-62){$- \vartheta V^{2,0}_2$}
\put(25,-25){$\vartheta$}
\end{picture}
}
\put(70,0){
\begin{picture}(1,1)
\put(0,-25){\line(1,0){70}}
\put(35,-25){\vector(1,0){5}}
\put(0,-25){\circle*{2}}
\put(70,-25){\circle*{2}}
\put(35,-20){$t$}
\put(-25,-20){$-t V_1^{1,0}$}
\put(71,-20){$+t V_1^{1,1}$}
\end{picture}
}
\end{picture}
\end{center}
\vspace{10ex}
\caption{Here are shown the orientation of $2$-cube and values of bound differentials $\vartheta V^{n,\varepsilon}_i  = \partial_i^{n,\varepsilon} \vartheta$.}
\end{figure} 
\par In figure 1, we have sketchy shown basic cocycles on the torus: if the dotted line is crossed any edge of the cube then the cocycle take $1$ on this edge, and this cocyle take $0$ on others edges.
\par Let us consider the $\smile$-product of basic cocycles. Since $\vartheta \circ V^{2,\varepsilon}_i =h_{\lambda_{\{i\}}^\varepsilon} \circ \vartheta $, we get  (see (6) and figure 2.)
$$(\alpha \smile \beta)(\vartheta) = \alpha \left(\vartheta V_1^{2,0} \right)\cdot \beta \left(\vartheta V_2^{2,1} \right) -  \alpha \left(\vartheta V_2^{2,0} \right)\cdot \beta \left(\vartheta V_1^{2,1} \right) = 0 \cdot 0 - (-1) \cdot (-1) = -1.$$
\par Thus, $\beta \smile \alpha$ --- is a basic cocycle of $H^2(\mathbb{T}^2; \mathbb{Z}).$ Further
$$(\beta \smile \alpha)(\vartheta) = \beta \left(\vartheta V_1^{2,0} \right)\cdot \alpha \left(\vartheta V_2^{2,1} \right) -  \beta \left(\vartheta V_2^{2,0} \right)\cdot \alpha \left(\vartheta V_1^{2,1} \right) = 1 \cdot 1 - 0 \cdot 0 = 1$$
\par So, we see that the cohomology ring $H^*(\mathbb{T}^2, \mathbb{Z})$ can be identified with the exterior algebra over the $\mathbb{Z}$-module $\mathbb{Z}$ whose generators are $\alpha$ and $\beta$.

\end{example}

\section*{Concluding Remark}
So, let us to sum up. For any precubical set $X \in \square_+^{op}\mathrm{Ens}$ and for any ring $R$ we get a graded cohomomology ring $H^*(X;R)$. If the ring $R$ has the unit then the ring $H^*(X;R)$ has the same unit. Further, if the ring $R$ is commutative then the ring $H^*(X;R)$ is anticommutative.

\end{document}